\documentclass{amsart}
\newtheorem{defin}{Definition}
\newtheorem{lemma}{Lemma}
\newtheorem{propo}{Proposition}
\newtheorem{theor}{Theorem}

\newcommand{\AAA}{A}
\newcommand{\ADM}[1]{A(#1)}

\newcommand{\BBB}{B}
\newcommand{\BOU}[1]{L(#1)}
\newcommand{\BSG}{B}

\newcommand{\CAR}{H}
\newcommand{\CHT}{\phi}

\newcommand{\COR}[1]{C(#1)}

\newcommand{\DEC}{\tilde{\MOS}}
\newcommand{\DET}[1]{\tilde{\MOS}(#1)}

\newcommand{\EDG}[2]{E_{#2}(#1)}
\newcommand{\EPS}{\epsilon}

\newcommand{\FAC}[1]{F(#1)}

\newcommand{\FRS}[1]{p(#1)}

\newcommand{\GRF}[1]{\Gamma(#1)}
\newcommand{\GRP}{G}

\newcommand{\HOM}{\mathrm{Hom}(\pi_1,\GRP)}

\newcommand{\IDA}{\mathcal{A}_\SFC}
\newcommand{\IDT}{\Delta_\SFC}
\newcommand{\IDS}{\mathcal{W}}
\newcommand{\INQ}{\psi}
\newcommand{\INT}{\mathbb{Z}}

\newcommand{\KAP}{\kappa}

\newcommand{\LOP}{\ell}

\newcommand{\MOS}{\mathcal{M}}
\newcommand{\MOT}[1]{\mathcal{M}(#1)}

\newcommand{\NAD}[1]{\bar{A}(#1)}
\newcommand{\NOR}[1]{N(#1)}

\newcommand{\OME}{\Omega}
\newcommand{\OMS}{\Omega_{\mathrm{s}}}

\newcommand{\OPP}{O}

\newcommand{\PAI}{I}
\newcommand{\PAT}[1]{C(#1)}
\newcommand{\PBL}{\mathcal{R}}

\newcommand{\PPP}[1]{P_{#1}}
\newcommand{\PSG}{T}

\newcommand{\QQQ}[1]{Q_{#1}}

\newcommand{\RLS}{\mathbb{R}}

\newcommand{\SCD}[1]{q(#1)}
\newcommand{\SFC}{\Sigma}
\newcommand{\SHT}{\hat{\Sigma}}
\newcommand{\SSS}{S}
\newcommand{\SST}{\Pi}

\newcommand{\TED}{\tilde{\mathcal{T}}}
\newcommand{\TEI}{\mathcal{T}}
\newcommand{\THE}{\theta}
\newcommand{\TRI}{T}
\newcommand{\TTT}{\rho}

\begin{document}

\title{Coordinates for the moduli space of  flat $PSL(2,\RLS)$-connections}

\author{R. M. Kashaev}
\thanks{This work is supported in part by the Swiss National Science
Foundation} 

\address{Section de math\'ematiques, Universit\'e de Gen\`eve
CP 240, 2-4, rue du Li\`evre, CH 1211 Gen\`eve 24, Suisse}

\address{V.A. Steklov Institute of Mathematics at St. Petersburg , 27, 
Fontanka, St. Petersburg 191023, Russia}
 
\email{Rinat.Kashaev@math.unige.ch}
\begin{abstract}
Let $\MOS$ be the moduli
space of irreducible flat $PSL(2,\RLS)$ connections on a 
punctured surface of finite type
with parabolic holonomies around punctures.
By using a notion of \emph{admissibility} of an ideal arc,
$\MOS$ is covered by dense open subsets associated to ideal
triangulations of the surface. A principal 
bundle over 
$\MOS$ is constructed which, when restricted to
the Teichm\"uller component of $\MOS$, is
isomorphic to the decorated
Teichm\"uller space of Penner. The construction gives a generalization to
$\MOS$ of Penner's coordinates for the Teichm\"uller space.
\end{abstract}

\maketitle

\section{Introduction}

In this paper we consider a punctured surface of
finite type
$\SFC=\SFC_{g,s}$ which is
the complement of a finite set
 of points $V\subset\overline{\SFC}$, $|V|=s$, called \emph{punctures},
in a closed oriented 
surface $\overline{\SFC}$ of genus $g$:
\[
\SFC=\overline{\SFC}\setminus V
\]
We
assume the following restrictions on $g$ and $s$:
\[
s>0,\quad \KAP\equiv 2g-2+s>0
\]
Denote  by $\GRP$ the group $PSL(2,\RLS)=SL(2,\RLS)/\{\pm1\}$. 
Let $\MOS$  be the moduli
space of irreducible flat $\GRP$-connections on $\SFC$ 
with parabolic holonomies around the punctures. Connections, representing
points of $\MOS$ will be called \emph{flat connections}.  There is a one
to one correspondence between elements of $\MOS$ 
and equivalence classes of irreducible
representations in $\GRP$ of the fundamental group of the surface 
with parabolicity conditions at
punctures. In other words, the moduli space
$\MOS$ can be identified with the quotient space $\HOM/\GRP$, where 
the group $\GRP$ acts  by conjugations freely and
properly so that $\MOS$ has the structure of a smooth manifold 
as the base of a
principal $\GRP$-bundle.  

In this paper we  extend the Penner's construction of the decorated
Teichm\"uller  space \cite{Penner1} to the whole moduli space $\MOS$.
Motivation for this work comes from nice properties of the decorated
Teichm\"uller space (such as simple and explicit realization of the
mapping class group action, simple form of the Weil--Petersson
symplectic structure, etc.) which are very useful for quantizing the
Teichm\"uller space \cite{CF,Kas}. 
The main results of this paper are formulated below in 
Theorems~\ref{t:1} and \ref{t:2}. 

With any path $a$ in $\SFC$, starting at some $x\in\SFC$ and ending at some
puncture $P$, we associate the homotopy class $\LOP(a)$ of loops
based at $x$ which go along $a$ towards $P$, then go around $P$ along
a small circle, and then return back to $x$ along $a^{-1}$.

Let $a$ and $b$ be two paths starting from one and the
  same point in 
$\SFC$ and ending at some punctures. Let $m\in\MOS$ be represented
by flat connection $h$. 
\begin{defin} 
The homotopy class of the path $a^{-1}b$ is called
  $m$-\emph{admissible} if 
$h$-holonomies along loops $\LOP(a)$ and $\LOP(b)$ 
belong to distinct parabolic
subgroups of $PSL(2,\RLS)$.
\end{defin}
It is clear that this notion depends on only the homotopy class
$[a^{-1}b]$ 
and moduli $m=[h]$, but not on
the choice of their representatives.

Recall that an \emph{ideal arc} on $\SFC$ is a nontrivial isotopy class of a
simple path running between punctures, and an \emph{ideal
  triangulation} of $\SFC$ is a maximal family of pairwise
nonintersecting ideal arcs. The sets $\IDA$ and $\IDT$ of ideal arcs
and ideal triangulations, respectively, are countable infinite.

With each ideal triangulation $\tau$ we can now associate an open dense
subset $\MOS_\tau$ of $\MOS$ defined as those  moduli of flat connections for
which $\tau$ is an admissible ideal triangulation, i.e. with all ideal
arcs admissible. 
\begin{theor}\label{t:1}
 Let $\SFC,\MOS,\{\MOS_{\tau}\}_{\tau\in\IDT}$ be as above.
 Then
\begin{description}
\item[(i)]the collection 
$\{\MOS_\tau\}_{\tau\in\IDT}$ 
is a covering for $\MOS$;
\item[(ii)] there exists a finite subcovering
$
\{\MOS_\tau\}_{\tau\in\SST}$,
$\SST\subset\IDT$, 
$|\SST|<\infty$.
\end{description}
\end{theor}

\begin{theor} \label{t:2}
Let $\SFC,\MOS,\{\MOS_\tau\}_{\tau\in\IDT}$ be as above.
 Then there exists a principal $\RLS_{>0}^s$-bundle 
$\pi\colon\DEC\rightarrow\MOS$ such that
\begin{description} 
\item[(i)] for each $\tau\in\IDT$ one has 
\(
\pi^{-1}(\MOS_\tau)=\coprod_{\EPS\in\{\pm1\}^{2\KAP}}\PBL_{\EPS}(\tau)
\)
with each $\PBL_{\EPS}(\tau)$ being a
principal $\RLS_{>0}^s$-bundle homeomorphic to the complement of
$\INQ_{\tau,\EPS}^{-1}(0)$ for certain rational mapping
$\INQ_{\tau,\EPS}\colon \RLS_{>0}^{3\KAP}\to\RLS$;

\item[(ii)] if $\sigma\colon \{\pm1\}^{2\KAP}\to \INT$ is defined
by 
\(
\sigma(\EPS)=\frac12\sum_{i=1}^{2\KAP}\EPS_i,
\)
then for $-\KAP\le k\le \KAP$, the sets
\(
\PBL_{k}=\bigcup_{\tau\in\IDT}\coprod_{\EPS\in\sigma^{-1}(k)}\PBL_{\EPS}(\tau)
\)
are principal $\RLS_{>0}^s$-bundles disjoint for different $k$;
\item[(iii)]  there exist
principal bundle isomorphisms
between $\PBL_{\pm\KAP}$ and the
decorated Teichm\"uller space $\TED$ of Penner. 
\end{description}
\end{theor}
 Theorem~\ref{t:1}(i) and Theorem~\ref{t:2}(ii)  imply that
\[
\MOS=\coprod_{-\KAP\le k\le \KAP}\MOS_k,\quad
\MOS_k=\PBL_k/\RLS_{>0}^s
\]
 The total number of these components is
$2\KAP+1$ which formally at $s=0$ coincides with the number of connected
components of the moduli space of flat $PSL(2,\RLS)$ connections on
a closed surface \cite{Gold}. In our case, however, the components $\MOS_k$ 
can be empty or 
not connected. For example, in the case $\KAP=1$,
$\MOS_{0}=\emptyset$. Part (iii) of Theorem~\ref{t:2} implies that
$\MOS_{\pm\KAP}$, being homeomorphic to the Teichm\"uller space, are
open cells of dimension $3\KAP-s$. 

Part (i) of Theorem~\ref{t:2} extends to $\MOS$ Penner's coordinate charts
for the Teichm\"uller space. Coordinatization is
given by assigning positive real numbers 
to edges and signs to faces of an ideal triangulation. Part (iii)
implies
that the
Teichm\"uller component of the moduli space corresponds to putting
all signs to one and the same value. In this
parameterization one still has an explicit description of both the action
of the mapping class group in $\MOS$ and the extension of the Weil--Petersson
symplectic form to $\MOS$. The latter is given
by the same formula as that of Penner in \cite{Penner2} 
for the Weil--Petersson form
(the pull-back thereof).

Restriction to parabolic holonomies around the punctures is essential
for our construction. Nevertheless, any other case, 
including the case of closed
surfaces, can also be considered in our framework by fixing holonomies
around some essential loops. The fact, that we are considering the whole
moduli space ensures that holonomies along essential loops
can be of any type: hyperbolic, parabolic or elliptic.

The rest of this paper is organized as follows. In Section~\ref{s:1}
we collect the material to be used in proving the theorems. Then,
Sections~\ref{s:2} and \ref{s:3} contain the proofs of two parts of
Theorem~\ref{t:1}. In the last Section~\ref{s:4} we prove Theorem~\ref{t:2}.
  
The author would like to thank A. Alekseev, A. Haefliger,
J.-C. Hausmann, and N. Reshetikhin for helpful discussions.  

\section{Preliminaries}\label{s:1}
\subsection{Subgroups of $PSL(2,\RLS)$}
Fix a parabolic subgroup $\PSG\simeq\RLS$ of $\GRP$. It's normalizer
$\BSG=\NOR{\PSG}$  is a Borel
subgroup 
in $\GRP$. Naturally, they are included in the exact sequence 
of group homomorphisms
\[
1\to\PSG\to\BSG\to\CAR\to 1
\]
where  $\CAR=\BSG/\PSG\simeq\RLS_{>0}$ is isomorphic to
an abelian subgroup of $\BSG$.
 The Bruhat
decomposition of $\GRP$ with respect to $\BSG$ consists of only two
cells
\[
\GRP=\BSG\THE \BSG\sqcup\BSG
\]
where $\THE\in\NOR{\CAR}\subset\GRP$ is a fixed representative 
of the only nontrivial element of the Weyl group $\NOR{\CAR}/\CAR$. 
We choose it in a unique
way by the condition $\THE^2=1$.

\subsection{Ideal triangulations}
 
For any $x,y\in\overline{\SFC}$ denote by $\OME(x,y)$ 
the set of paths in $\SFC$ starting
at $x$ and ending
at $y$. The subset $\OMS(x,y)\subset\OME(x,y)$ consists of simple paths. 
Denote also $\OME(x)=\OME(x,x)$ and $\OMS(x)=\OMS(x,x)$. Two paths
are called \emph{not intersecting} if their interior parts 
do not intersect. 

\emph{Ideal arcs} are
nontrivial isotopy classes in $\OMS(P,Q)$, for all $P,Q\in V$.  
Associated to each ideal arc $e$ there are two 
vertices $\PPP{e},\QQQ{e}\in V$, possibly coinciding, such that $e$ is
an isotopy class in $\OMS(\PPP{e},\QQQ{e})$. The set of all ideal
arcs is denoted $\IDA$. 

An \emph{ideal
  triangulation} of $\SFC$ is a maximal set of pairwise nonintersecting ideal
arcs.  The  set of all ideal triangulations of $\SFC$ is denoted by $\IDT$.

Let an ideal triangulation be represented by a system of simple
  pairwise nonintersecting curves. It is clear that the complement of
  $\SFC$ 
to such  system is a disjoint union of cells
  homeomorphic to triangles. The isotopy classes of such triangles will
  be called \emph{ideal triangles}.
Each ideal triangulation $\tau$ is thus associated with a cell complex
homeomorphic to $\overline{\SFC}$ with the set of
punctures  $V$ as vertices, the ideal arcs $\EDG{\tau}{}=\tau$ as (oriented)
edges, and
ideal triangles $\FAC{\tau}=\bar\tau$ as faces. In addition, we denote by 
$\PAT{\tau}=\coprod_{P\in V}\PAT{\tau,P}$ the set of oriented 
corners of ideal triangles, i.e. the connected components of 
intersections of ideal 
triangles  with small disks centered at punctures, where 
$\PAT{\tau,P}$ is the set of corners containing
vertex $P$. In fact, one has $s$ vertices,
$3\KAP$ edges, $2\KAP$ faces, and $6\KAP$ corners in any ideal
triangulation of $\SFC$.
Each corner $c\in\PAT{\tau,P}$ belongs to unique triangle
$\TRI_c$ with the side $\OPP_c\in \tau$ opposite to $P$.  
Denote two other sides of $\TRI_c$ by $\AAA_c$ and $\BBB_c$. 
For an ideal polygon $p$ in $\SFC$ we denote by $\EDG{p}{}$ and
  $\COR{p}$ the sets of it's sides and corners, respectively.

If $a\in\OME(x,P)$, where $x\in\overline{\SFC}$ and $P\in V$, 
we denote by $\LOP(a)$ the homotopy class of the loop based at $x$ and
obtained by going along $a$ towards $P$, going around $P$ along
 a small circle in the direction induced from the orientation of
 $\SFC$, and returning back to $x$ along $a^{-1}$. In the case, when
 $x=Q\in V$, the class $\LOP(a)$, when it corresponds
 to a simple loop, is naturally associated
 with an ideal arc based at $Q$.
\subsection{Graph connections}
For $\alpha\subset\IDA$ define an open subset $\MOT{\alpha}\subset\MOS$ by the
condition that for any $m\in \MOT{\alpha}$ the set $\alpha$ is the
maximal family of $m$-admissible ideal arcs. To avoid trivial cases with
$\MOT{\alpha}=\emptyset$, we define also
\[
\IDS=\left\{\alpha\subset\IDA\left|\,\MOT{\alpha}\ne\emptyset\right.\right\}
\] 
It is clear that
\[
\MOS=\coprod_{\alpha\in\IDS}\MOT{\alpha}
\]
and the Teichm\"uller component $\TEI\subset\MOT{\IDA}$. 
Theorem~\ref{t:1} partially characterizes
elements of
$\IDS$, namely each of them contains at least one ideal
triangulation. It is clear also that
\[
\MOS_\tau=\coprod_{\alpha\in\IDS(\tau)}\MOT{\alpha},\quad
\IDS(\tau)=\{\alpha\in\IDS|\,\tau\subset\alpha\}
\]

Define a pairing
\[
\PAI\colon\IDA\times\IDA\to\INT_{\ge0},\quad
\PAI(e,f)=\min(|a\cap b|\;a\in e,b\in f)
\]
Let $\SHT\subset\SFC$ be the complement of a disjoint union of small
open disks centered at the punctures. The boundary of $\SHT$ is a
disjoint union of $s$ circles $\{\BOU{P}\}_{P\in V}$
associated with punctures
\[
\partial\SHT=\coprod_{P\in V}\BOU{P}
\]
For each $\alpha\in\IDS$ we associate a set
$\GRF{\alpha}=\{(e,\FRS{e},\SCD{e})\}_{e\in\alpha}$, where $\FRS{e}\in
\BOU{\PPP{e}}$ and $\SCD{e}\in
\BOU{\QQQ{e}}$ are chosen so that for any ideal arcs $e$ and
$f$ there exist such $a\in e$ and
$b\in f$, that $a\cap\SHT\in\OMS(\FRS{e},\SCD{e})$,
 $b\cap\SHT\in\OMS(\FRS{f},\SCD{f})$, and $|a\cap
 b|=\PAI(e,f)$. We shall think of $\GRF{\alpha}$ as a graph in $\SFC$ with
 points
 $\{\FRS{e},\SCD{e}\}_{e\in\alpha}$ as vertices and two
 types of edges: ``long'' edges corresponding to ideal arcs, and
``short'' edges corresponding to segments of the boundary
 components of $\SHT$ between vertices. If $\alpha\subset\beta$, we
 identify $\GRF{\alpha}$ with the corresponding subgraph of
 $\GRF{\beta}$. Thus, for any $\alpha\in\IDS$, $\GRF{\alpha}$ is a
 subgraph
of $\GRF{\IDA}$. 

Any polygon $p$ in $\GRF{\IDA}$ is uniquely associated with an ideal
polygon $p'$ in $\SFC$, the long and short sides of $p$  corresponding
to sides and corners of $p'$, respectively.

\begin{defin}
 A \emph{graph $\GRP$-connection} on  a graph $\Gamma\subset\SFC$ is an
 isomorphism class of $\GRP$-connections on $\SFC$ with fixed parallel
 transport operators along edges of $\Gamma$. A \emph{graph gauge
 transformation} of a graph connection is a usual gauge
 transformation  modulo gauge transformations relating different
 representatives of graph connections. 
\end{defin}
Equivalently, one can define a graph $\GRP$-connection as a
representation in $\GRP$ of the path-edge groupoid of $\Gamma$, graph gauge
 transformations being equivalence transformations of representations.
It is clear that if graph $\Gamma$ contains a subgraph
homotopically equivalent to $\SFC$, then there is a one to one
correspondence between equivalence classes of flat graph connections
on $\Gamma$ (with respect to graph gauge transformations) and moduli
of flat connections on $\SFC$, see, for example, \cite{fock_rosly}.
A typical example of a graph connection is given by a representation
of the fundamental group $\pi_1(\SFC,x)$ in $\GRP$, where the graph
$\Gamma$ is given by homotopy classes in $\OME(x)$ as edges and a
single vertex $x$. Graph gauge transformations in this case correspond
to overall conjugations by elements of $\GRP$.

\section{Proof of Theorem~\ref{t:1}(i)}\label{s:2}

\subsection{Preparation}
We shall find it useful the following lemmas, where 
in each of them it is assumed that a flat connection is given.
\begin{lemma}\label{l:1}
There exists an admissible ideal arc
 starting and ending at any given puncture.
\end{lemma}
\begin{proof}
Suppose that there are no such admissible ideal arcs. Fix a base point
$x\in \SFC$ and choose a puncture $P$ together with a curve
$b\in\OMS(x,P)$. Denote by $h_0$ 
the associated to $\LOP(b)$ holonomy.
Let $\delta\in\OMS(x)$ be such that it does not intersect $b$. Such 
loops generate the fundamental group $\pi_1(\SFC,x)$.
We associate to $\delta$ the loop $b^{-1}\delta b\in\OME(P)$, 
the homotopy class of which uniquely defines an ideal arc $e$
based at $P$. Choose a representative $a^{-1}b$ for $e$, where
$a\in\OME(x,P)$. Then the holonomies
corresponding to loops $\LOP(b)$ and $\LOP(a)$ 
are given by elements $h_0$ and $h_\delta^{-1} h_0
h_\delta$, respectively, 
where $h_\delta$ is the
holonomy associated to $\delta$.  The nonadmissibility
of $e$ implies that the elements $h_0$ and
$h_\delta^{-1}h_0 h_\delta$  belong to one and the same parabolic
subgroup, 
and thus 
the element $h_\delta$ belongs to this subgroup's normalizer. 
Conjugating the element $h_0$
to an upper triangular form, one simultaneously brings to an upper
triangular form also element
$h_\delta$. Taking into account the fact that
element $h_0$ remains the same
independently
of the choice  of the loop $\delta$, we conclude that
the
corresponding to the flat connection representation of the fundamental group
$\pi_1(\SFC,x)$ is brought to an upper triangular form. 
This contradicts the irreducibility of the flat connection.
\end{proof}
\begin{lemma}\label{l:2}
Let $P,Q,R\in V$ (possibly coinciding) and 
$a\in\OME(P,Q)$, $b\in\OME(Q,R)$, $c\in\OME(R,P)$ be such that $abc$
is homotopic to a trivial loop. Let the homotopy 
class $[a]$ be admissible. Then at least one of $[b]$ and $[c]$ is
also admissible. In particular,   
an ideal triangle, having at least one admissible side, has at least two
 admissible sides.
\end{lemma}
\begin{proof}
Let $x\in\SFC$.
 Choose $p\in\OME(x,P)$,
$q\in\OME(x,Q)$, $r\in\OME(x,R)$ so that
$a=p^{-1}q$, $b=q^{-1}r$, $c=r^{-1}p$. Admissibility of
$[a]$ means that the holonomies $h_{\LOP(p)}$ and
$h_{\LOP(q)}$ belong to distinct parabolic
subgroups of $\GRP$. Then the holonomy $h_{\LOP(r)}$ can
belong to only one of these subgroups.
\end{proof}
\begin{lemma}\label{l:4}
For any two distinct punctures there exists an
admissible ideal arc connecting them.
\end{lemma}
\begin{proof}
Let $P$ and $Q$ be two distinct punctures, and let $e$ be an
admissible 
ideal arc based at $P$ (which exists by Lemma~\ref{l:1}). 
Fix $x\in \SFC$ and $c\in\OMS(x,Q)$ so that  
$e$ is represented by $a^{-1}b$, where $a,b\in\OMS(x,P)$ do not
intersect $c$. The homotopy classes $[a^{-1}c]$ and $[b^{-1}c]$  define
two (possibly coinciding) ideal arcs, connecting $P$ and $Q$. By 
Lemma~\ref{l:2} (applied to $e$, $[b^{-1}c]$ and $[c^{-1}a]$), 
at least one of them is admissible. 
\end{proof}
\begin{lemma}\label{l:5}
Let $P,Q\in V$ (possibly coinciding) and $d\in\OME(P,Q)$. Then $[d]$ is
admissible if and only if
$\LOP(d)$ is admissible.  
\end{lemma}
\begin{proof}
 Let $x\in\SFC$, $a,b\in\OME(x,P)$, $c\in\OME(x,Q)$
be such
that $\LOP(d)=[a^{-1}b]$ and $[d]=[a^{-1}c]=[b^{-1}c]$. 
We have the identity
\begin{equation}\label{eq:ll}
\LOP(c)\LOP(a)=\LOP(b)\LOP(c)=\LOP(cd^{-1})
\end{equation} 
Suppose $\LOP(d)$ is nonadmissible. Then holonomies
$h_{\LOP(a)}$ and $h_{\LOP(b)}$  belong to one and the
same parabolic subgroup $T\subset\GRP$.
Then eqn~\eqref{eq:ll} implies
that $h_{\LOP(c)}$ is in 
the normalizer of $T$, so 
$h_{\LOP(c)}\in T$ since it is itself a parabolic element. Thus,
$[d]$ is nonadmissible. Conversely, if $[d]$ is
nonadmissible, then the holonomies around $\LOP(a)$, $\LOP(b)$ and 
$\LOP(c)$ are in
one and the same parabolic subgroup, i.e. $\LOP(d)$ is nonadmissible too.
\end{proof}
\begin{lemma}\label{l:8}
On $\SFC$ there exists a tree with the set of punctures
 as vertices and $s-1$ pairwise
nonintersecting admissible ideal arcs as edges.
\end{lemma}
\begin{proof}
 Enumerate the punctures
    $P_1,\ldots,P_s$. 
By Lemma~\ref{l:4}, there exists an admissible ideal arc connecting $P_1$ and
    $P_2$. Denote  by $\Gamma_2$ the tree
with two vertices $P_1$ and $P_2$ and the
    admissible ideal arc as it's edge. 
Assume that the punctures $P_1,\ldots,P_k$ are
    vertices of a tree $\Gamma_{k}$ with $k-1$ edges being pairwise
    nonintersecting admissible ideal arcs, where $1<k<s$. It is clear that
one can find 
such an  edge $e$ of $\Gamma_{k}$, connecting vertices $P_i$ and
    $P_j$, that there
    exists
an ideal triangle $t$ with vertices $P_i,P_j,P_{k+1}$ and whose edges
    do not
    intersect the edges of $\Gamma_{k}$. By Lemma~\ref{l:2}
    one of the admissible edges of $t$  connects $P_{k+1}$
either to $P_i$ or $P_j$. By adding this edge and $P_{k+1}$ to
    $\Gamma_{k}$ we obtain a tree $\Gamma_{k+1}$ with vertices
$P_1,\ldots,P_{k+1}$ and $k$ pairwise
    nonintersecting admissible ideal arcs as edges. By recurrence,
    $\Gamma_s$ is a tree with the required properties. 
\end{proof}
\begin{defin}
A \emph{short diagonal} of an ideal $n$-gon $p$ is an ideal arc in $p$
whose complement  in  $p$ is an ideal triangle and an ideal $(n-1)$-gon.
\end{defin}
\begin{lemma}\label{l:3}
There exists such a system
of $\KAP+1$ pairwise nonintersecting  admissible ideal arcs,
containing all punctures, that
the complement in $\SFC$ to this system is an ideal 
$(2\KAP+2)$-gon. For any such system the corresponding polygon
 has at least one admissible short diagonal.
\end{lemma}
\begin{proof}
Let $\Gamma_s$ be the tree of Lemma~\ref{l:8}, and let $P$ be a puncture.
There exists such a system $\Gamma'$ of $2g$ pairwise nonintersecting
 ideal arcs (not necessarily
admissible), based
at  $P$ and not intersecting the edges of tree $\Gamma_s$, that 
the  complement of $\Gamma'$ in $\SFC$ is an $4g$-gon with $\Gamma_s$
inside of it and attached to one of it's vertices. This means that the
complement of the graph
$\Gamma^0=\Gamma_s\cup\Gamma'$ is an ideal polygon $p_0$ with
$(2\KAP+2)$ sides.
In fact, $\Gamma^0$ contains at least one admissible ideal arc. This is
clear for $s>1$ as the $s-1$ edges of the subgraph $\Gamma_s$ are
admissible. In the case, when $s=1$, one can choose $\Gamma'$ with 
at least one admissible loop (this is possible due to
Lemma~\ref{l:1}). Thus, we can assume that the polygon $p_0$ 
has nonempty set of admissible sides.
Suppose $p_0$ has $2n$
nonadmissible sides, where  $1\le n<\KAP+1$. 
Let $a$ be an admissible side next to a
nonadmissible side $b$. Let $c$ be the short
diagonal of $p_0$ which together with $a$ and $b$ bounds an ideal 
triangle. By Lemma~\ref{l:2}, $c$ is admissible. Replacing $b$ by $c$
in $\Gamma^0$, we
obtain another graph $\Gamma^1$ corresponding to ideal polygon $p_1$ 
with $2n-2$ nonadmissible sides. By repeating this procedure recurrently
$n-1$ more times we obtain a graph $\Gamma^n$ and the corresponding
polygon $p_n$
whose all $2\KAP+2$ sides are admissible. 

Let now  $\Gamma$ be any such graph with the corresponding polygon
$p$. Assume that $p$
does not have admissible short diagonals. Let
$P\in V$ and $e_0,e_1,\ldots,e_{k-1}$ be
the edges of $\Gamma$ (possibly with repetitions) which intersect
a small circle around $P$ with the cyclic order
induced from the orientation of $\SFC$. Let $e_i=[d_i]$ for some
$d_i\in\OME(P_i,P)$, $0\le i<k$. Fix $x\in\SFC$ and $a_i\in\OME(x,P_i)$,
$0\le i\le k$, $P_k=P_0$, so that
$\{ a_i^{-1}a_{i+1}\}_{0\le i<k}$ represent  $k$ 
shorts diagonals (possibly with repetitions) of $p$.
Since all short
diagonals of $p$ are nonadmissible, then
the  holonomies around 
$\LOP(a_i)$, $0\le i\le
k$, belong to one and the same parabolic subgroup of $\GRP$. This
contradicts  Lemma~\ref{l:5} according to which
$\LOP(d_0)=[a_0^{-1}a_k]$ is admissible.
 \end{proof}
\begin{lemma}\label{l:6}
An ideal $n$-gon $p$, where $n>4$, with all 
sides admissible and at least one admissible short diagonal
has such a pair of nonintersecting admissible diagonals, that it's complement
in $p$ is a union of an ideal $(n-2)$-gon and two ideal triangles.
\end{lemma}
\begin{proof}
As a pair of  admissible nonintersecting short diagonals has the 
required property, we
assume that there are no such pairs. Let $d$ be an admissible short
diagonal with $p\setminus d=q\sqcup t$, where $q$ is an ideal
$(n-1)$-gon and $t$ an ideal triangle. Then $q$ has at least
one admissible short diagonal, since otherwise, by Lemma~\ref{l:2}, 
the two nonintersecting short
diagonals of $p$ containing the puncture opposite to $d$ in $t$ are
admissible.
\end{proof}

\begin{lemma}\label{l:7}
An ideal $n$-gon, where $n>3$, with all 
sides admissible and at least one admissible short diagonal,
has $n-3$ pairwise nonintersecting admissible diagonals.
\end{lemma}
\begin{proof}
The case $n=4$ is automatic.
The case $n=5$ is equivalent to Lemma~\ref{l:6}. Assume that
$n=5+k$, where $k>0$.
Denote the initial $n$-gon by $p_0$. Of the 
two diagonals of Lemma~\ref{l:6}, at least one is short. Let it
be $d_1$ in $p_0$. Then $p_0\setminus d_1=p_1\sqcup t_1$, where $t_1$ is an
ideal triangle and $p_1$ an ideal $(n-1)$-gon. Polygon $p_1$ itself
satisfies the 
conditions of Lemma~\ref{l:6}. By recurrence, we obtain a sequence of pairwise
nonintersecting admissible diagonals $\{ d_i\}_{1\le i\le k}$, ideal
triangles $\{ t_i\}_{1\le i\le k}$, and a pentagon $p_k$, satisfying the
conditions
of Lemma~\ref{l:6}.
\end{proof}

\subsection{Proof of part (i)}
This part of the theorem is equivalent to the following
\begin{propo}\label{p:1}
For each flat connection there exists an admissible ideal 
triangulation.
\end{propo}
\begin{proof}
Given a flat connection. Let $\Gamma$ be the system of $\KAP+1$
ideal arcs of Lemma~\ref{l:3}. By Lemma~\ref{l:7} we can complete
$\Gamma$ to a system
$\tilde{\Gamma}$ of $\KAP+1+2(\KAP+1)-3=3\KAP$ 
pairwise nonintersecting admissible ideal arcs.  $\tilde{\Gamma}$
is evidently an admissible ideal triangulation. 
\end{proof}

\section{Proof of Theorem~\ref{t:1}(ii)}\label{s:3}
\subsection{Preparation}
Let $\tau\in\IDT$. Consider
an ideal quadrilateral $q$ in $\tau$ 
with diagonal $e\in \tau$. Replacing $e$ by
another diagonal $e'$ of $q$ (such operation is called a \emph{flip on
$e$}), 
we obtain another ideal triangulation
$\tau'$. 

In the rest of this subsection we assume that
a flat connection is given.
\begin{lemma}\label{ll:1}
Any ideal triangulation of $\SFC$ has at least one admissible edge. 
\end{lemma}
\begin{proof}
Let $\tau\in\IDT$. Assume that all edges of $\tau$ are
nonadmissible. Then, due to Lemma~\ref{l:2}, for any ideal
quadrilateral 
in $\tau$ both it's
diagonals are nonadmissible. Thus, any flip in $\tau$ leads to an ideal
triangulation with all edges nonadmissible. As any two ideal
triangulations can be related to each other by a finite sequence of
flips \cite{Harer,Mosher}, there are no
admissible ideal triangulations on $\SFC$. This contradicts
Theorem~\ref{t:1}(i). 
\end{proof}
\begin{defin}
An \emph{almost admissible} ideal triangle is an ideal triangle
with at most one nonadmissible side.
\end{defin}

\begin{lemma}\label{ll:2}
Any ideal triangulation of $\SFC$ can be replaced by an ideal
triangulation with all ideal triangles almost admissible.
\end{lemma}

\begin{proof}
Let $\tau_0\in\IDT$,
$\ADM{\tau_0}$ be the subset of 
almost admissible ideal triangles, and
$\NAD{\tau_0}=\FAC{\tau_0}\setminus\ADM{\tau_0}$.
By Lemma~\ref{ll:1}, $\ADM{\tau_0}\ne\emptyset$.
Assume $k=|\NAD{\tau_0}|>0$.   
Let $e\in \tau_0$ be such that the quadrilateral $q$, containing
$e$ as a diagonal, is composed of triangles
$t_1\in \NAD{\tau_0}$ and $t_2\in\ADM{\tau_0}$.  Then, 
by Lemma~\ref{l:2}, $t_2$ has exactly two admissible edges, and the second 
 diagonal $e'$ of $q$ is admissible. Replacing $e$ by $e'$, we
obtain ideal triangulation $\tau_1$ with
$|\NAD{\tau_1}|=k-1$. 
By induction, for ideal triangulation $\tau_k$ one has
$\NAD{\tau_k}=\emptyset$. 
\end{proof}

\subsection{Proof of part (ii)}

Fix $\tau\in\IDT$.
For any $f\colon\tau\to\{\pm1\}$
associate a subset 
$\MOS(f)\subset\MOS$ given by collection of those $m\in\MOS$, for which
edge $e\in \tau$ is admissible if $f(e)=1$ and nonadmissible if
$f(e)=-1$. It is clear that
\[
\MOS=\coprod_{f\in\SSS(\tau)}\MOS(f)
\]
where $\SSS(\tau)$ is the collection of those
$f$ for which $\MOS(f)\ne\emptyset$. Evidently,
$|\SSS(\tau)|< 2^{3\KAP}$. Now,  Theorem~\ref{t:1}(ii)
directly follows from
\begin{propo}
For each $f\in\SSS(\tau)$ one has a finite collection of
ideal triangulations 
\[
\TTT(f)\subset\IDT,\quad 
|\TTT(f)|=\frac1{2\KAP+1}\binom{4\KAP}{2\KAP}
\]
such that
\(
\MOS(f)\subset\bigcup_{\tau\in\TTT(f)}\MOS_{\tau}.
\)
\end{propo}
\begin{proof}
Given $f \in\SSS(\tau)$.
By using Lemma~\ref{ll:2}, we replace $\tau$ by $\tau_f\in\IDT$
such that for each $m\in\MOS(f)$ all faces of $\tau_f$ are almost admissible.
It is easily seen that there exists such a subset $A\subset \tau_f$
of $\KAP+1$  edges that the complement of $A$ in $\SFC$ is an ideal
$(2\KAP+2)$-gon  $p_A$ with all sides admissible for any $m\in\MOS(f)$.
To each ideal triangulation of $p_A$ there corresponds a unique ideal
triangulation of $\SFC$ with $A$ as a subset  of ideal arcs. 
We define $\TTT(f)$ as the collection of all such ideal
triangulations. The size of the set $\TTT(f)$ is equal to the
number of triangulations of a $(2\KAP+2)$-gon, i.e. the Catalan number
\[
\frac1{2\KAP+1}\binom{4\KAP}{2\KAP}
\]
Let $m\in\MOS(f)$. By construction, all sides of $p_A$ are
$m$-admissible, and, according to Lemma~\ref{l:3}, $p_A$ has at least
one $m$-admissible short diagonal. Then, by Lemma~\ref{l:7}, there exists 
an $m$-admissible triangulation of $p_A$, i.e. there exists such
$\tau\in\TTT(f)$ that $m\in\MOS_\tau$.
\end{proof}
It remains only to note that in Theorem~\ref{t:1}(ii) one can take
$\Pi=\bigcup_{f\in\SSS(\tau)}\TTT(f)$. It is clear that 
\[
|\Pi|<\frac{2^{3\KAP}}{2\KAP+1}\binom{4\KAP}{2\KAP}
\]
This is very rough estimation. In fact, for example, 
in the case $\SFC=\SFC_{0,4}$ one can show that there exists
$\Pi_{\mathrm{min}}$ of size $|\Pi_{\mathrm{min}}|=7$.

\section{Proof of Theorem~\ref{t:2}}\label{s:4}
We construct a fiber bundle
$\pi\colon\DEC\to\MOS$ as the disjoint
union of fiber bundles
$\pi\colon\DET{\alpha}\to\MOT{\alpha}$ for all $\alpha\in\IDS$. The
 fiber $\pi^{-1}(m)$ over a point $m\in\MOT{\alpha}$
consists of flat graph connections on $\GRF{\alpha}$ which represent $m$ and
whose parallel transport operators  belong
to the parabolic subgroup $\PSG$ for all short edges, and to the subset 
$\THE\CAR$ for all long edges.
Let us show that $\DEC$ is in fact a principal $\CAR^s$-bundle over $\MOS$.

Given $\alpha\in\IDS$ and $m\in\MOT{\alpha}$. Our strategy is to
show that any flat connection $h$, representing $m$, is equivalent to any 
flat graph connection in $\pi^{-1}(m)$.

 As every vertex  of $\GRF{\alpha}$ belongs to a
boundary component of $\SHT$, the
$h$-holonomies around boundary components, based at vertices of
$\GRF{\alpha}$, are parabolic.
 Thus, one always can replace $h$ by such an equivalent connection that all
associated holonomies around the boundary components, based at vertices,
belong to one and the same parabolic subgroup $\PSG\subset\GRP$. This
automatically 
makes the parallel transport operators along short edges
to belong to $\BSG=\NOR{\PSG}$.
The parallel transport operators along
long edges belong to the large Bruhat cell $\BSG\THE\BSG$. Thus, 
making an appropriate graph gauge transformation with values in $\BSG$, we
can make  the parallel transport operators along  all long edges to
belong to the subset 
$\THE\CAR$. The remaining freedom in graph gauge transformations
consists of arbitrary $\CAR$ valued functions on vertices of
$\GRF{\alpha}$, 
which can be used
to eliminate the $\CAR$-parts of parallel transport operators 
along all short edges, thus obtaining  elements of $\pi^{-1}(m)$.
This still leaves unfixed the graph gauge transformations,
given by $\CAR$ valued functions taking one and the same value on
all vertices associated with one boundary component. This
residual gauge group is isomorphic to $\CAR^{s}$, 
which freely and transitively acts in the space
$\pi^{-1}(m)$. 
Thus, the union $\DET{\alpha}=\coprod_{m\in\MOT{\alpha}}\pi^{-1}(m)$
is a principal $\CAR^s$-bundle over $\MOT{\alpha}$.

\subsection{Proof of part (i)} 
\begin{lemma}\label{l:11}
Given $\alpha\in\IDS$ and $m\in\MOT{\alpha}$. Let $f$ be a hexagonal
face in $\GRF{\alpha}$ corresponding to ideal triangle $t$ in
$\SFC$. Then, 
there exists $\epsilon=\pm1$ such that the parallel
transport operators along short edges of $f$ are uniquely defined in
terms of $\epsilon$ and the parallel
transport operators along long edges of $f$.
\end{lemma}
\begin{proof}
Let us identify $\BSG$ with  the set $\RLS_{>0}\times\RLS$ where the 
Lie group structure is that of upper triangular two-by-two matrices of
unit determinant 
\[
(u_1,v_1)(u_2,v_2)=(u_1u_2,u_1v_2+v_1u_2^{-1}),\quad (u_i,v_i)\in
\RLS_{>0}\times\RLS,\quad i=1,2
\]
Then, the subgroups $\CAR$ and $\PSG$ are identified with the
components $\RLS_{>0}$ and $\RLS$ respectively. In this way the
parallel transport operators are associated with mapping
$u\colon\COR{t}\to\RLS$
for short edges, and $v\colon\EDG{t}{}\to\RLS_{>0}$ for long edges.
 Now, explicit calculation
of the condition that the holonomy around $\partial f$ is trivial
shows that the real number $u(c)$, associated with corner $c$ of $t$, is
given by the formula
\begin{equation}\label{eq:uv}
u(c)=\epsilon\frac{v(\OPP_c)}{v(\AAA_c)v(\BBB_c)}
\end{equation}
\end{proof}
\begin{lemma}\label{l:12}
Given $\alpha\in\IDS$ and $m\in\MOT{\alpha}$. Let $f$ be a hexagonal
face in $\GRF{\alpha}$ corresponding to an ideal triangle $t$ in
$\SFC$. Then, the parallel transport
operators along edges of $f$ are uniquely
determined in terms of  parallel transport
operators along any two
long edges and any one short edge of $f$. 
\end{lemma}
\begin{proof}
This is a direct consequence of Lemma~\ref{l:11}. Indeed, let $a,b,c$ 
be the corners of $t$, $a$ and $b$ being opposite to $\AAA_c$ and $\BBB_c$,
respectively.  Solving
eqns~\eqref{eq:uv}, for example, with respect to $v(\OPP_c)$, $u(a)$,
and $u(b)$, 
we obtain explicitly
\[
v(\OPP_c)=|u(c)|v(\AAA_c)v(\BBB_c),\quad
u(a)=\frac1{u(c)(v(\BBB_c))^2},\quad u(b)=\frac1{u(c)(v(\AAA_c))^2}
\]
It is clear that instead of $u(c)$ one can take  either $u(a)$ or $u(b)$.
\end{proof}
Let $\GRF{\alpha,\tau}$ be the smallest
subgraph of $\GRF{\alpha}$ containing all long edges 
corresponding to $\tau\in\IDT$.
The parallel transport operators along edges of
$\GRF{\alpha,\tau}$ uniquely determine all parallel transport operators along
edges of $\GRF{\alpha}$. This is a consequence of
Lemma~\ref{l:12}. Indeed, for any $e\in\alpha\setminus\tau$ there
exists such a finite sequence of ideal arcs
$\{e_i\}_{1\le i\le n}\subset\alpha$ with $e_n=e$,  that $e_i$ is
contained in an ideal triangle whose other two sides are contained in
$\tau\cup\{e_j\}_{1\le j<i}$. Using recurrently
Lemma~\ref{l:12}, one can thus calculate the parallel transport operators
along all edges of $\GRF{\alpha}$ in terms of those of $\GRF{\alpha,\tau}$.
Besides, by Lemma~\ref{l:11},  each
point of $\pi^{-1}(\MOS_\tau)$ is uniquely associated with
 an element of the set
$\RLS_{>0}^{\tau}\times\{\pm1\}^{\bar{\tau}}$, where we denote
by $X^Y$ the set of mappings from set $Y$ to set $X$,
$\tau=\EDG{\tau}{}$, $\bar{\tau}=\FAC{\tau}$.
\begin{propo}\label{p:3}
Let $\tau\in\IDT$ and mapping
\(
\CHT_\tau\colon
\RLS_{>0}^{\tau}\times\{\pm1\}^{\bar{\tau}}\to \RLS
\)
be defined by the formula
\[
\CHT_\tau(f,\EPS)=
\prod_{P\in V}\CHT_{\tau,P}(f,\EPS),\quad
\CHT_{\tau,P}(f,\EPS)=\sum_{c\in\PAT{\tau,P}}\EPS(\TRI_c)\frac{
f(\OPP_c)}{f(\AAA_c)f(\BBB_c)}, 
\]
Then, there is a one to one correspondence between $\pi^{-1}(\MOS_\tau)$ and
the complement of $\CHT_\tau^{-1}(0)$, 
the action of the structure group $\CAR^s$ being
realized by the action of the group $\RLS_{>0}^V\simeq\CAR^s$:
\[
\RLS_{>0}^{\tau}\times\{\pm1\}^{\bar{\tau}}\times\RLS_{>0}^V\to
\RLS_{>0}^{\tau}\times\{\pm1\}^{\bar{\tau}}
\]
\begin{equation}\label{eq:sg}
(f,\EPS,h)\mapsto (f^h,\EPS),\quad
f^h(e)=f(e)h(\PPP{e})h(\QQQ{e})
\end{equation}
\end{propo}  
\begin{proof}
One only needs to satisfy the $s$ conditions for the holonomies around
the punctures to be parabolic, which, by using Lemma~\ref{l:11}, take the form
$\CHT_{\tau}(f,\EPS)\ne0$.

Verification of the action of the structure group is straightforward.
\end{proof}
To complete the proof of part (i), it remains to notice that the 
component $\EPS$ in formula~\eqref{eq:sg} is not affected under the
action of the structure group, so that for each fixed $\EPS$
there corresponds
a subbundle $\PBL_{\EPS}(\tau)\subset\pi^{-1}(\MOS_\tau)$.
\subsection{Proof of part (ii)}
\begin{propo}\label{p:4}
Consider two ideal triangulations $\tau$ and $\tau'$ related by a
single flip on $e\in\tau$ so that ideal triangles
$t_1,t_2\in\bar{\tau}$ with sides 
$e,a,b$ and $e,c,d$, respectively (the cyclic order of the edges being induced
from the orientation of $\SFC$), are replaced by
ideal triangles $t_1',t_2'\in\bar{\tau}'$
with sides
$e',d,a$ and $e',b,c$, respectively. Then
two pairs
\[
(f,\EPS)\in\RLS_{>0}^{\tau}\times\{\pm1\}^{\bar{\tau}}\setminus
\CHT_{\tau}^{-1}(0),\quad
(f',\EPS')\in
\RLS_{>0}^{\tau'}\times\{\pm1\}^{\bar{\tau}'}\setminus
\CHT_{\tau'}^{-1}(0)
\] 
correspond to 
one and the same element in $\pi^{-1}(\MOS_\tau\cap\MOS_{\tau'})$ if and
only if 
\begin{equation}\label{eq:fx}
f'(x)=f(x)\quad\forall\, x\ne e'
\end{equation}
\begin{equation}\label{eq:et}
\EPS'(t)=\EPS(t)\quad\forall\, t\ne t_i',\ i=1,2
\end{equation}
\begin{equation}\label{eq:ep}
\EPS'(t_1')\EPS'(t_2')=\EPS(t_1)\EPS(t_2)
\end{equation}
\begin{equation}\label{eq:ee}
\EPS'(t_1')f(e)f'(e')=\EPS(t_2)f(a)f(c)+\EPS(t_1)f(b)f(d)
\end{equation}

In particular,  if pair $(f,\EPS)$, representing
$m\in\MOS_\tau$, is such that
\[
\frac{f(a)f(c)}{f(b)f(d)}+\EPS(t_1)\EPS(t_2)=0
\] then
$m\not\in\MOS_{\tau'}$.
\end{propo}
\begin{proof}
This is a straightforward calculation by using
Lemma~\ref{l:11}. Indeed, let $\alpha\in\IDS$ be such that
$\tau\cup\tau'\subset\alpha$. 
Then, comparing the parallel transport operators along edges
of $\GRF{\alpha,\tau}\cap\GRF{\alpha,\tau'}$,
we immediately obtain the eqns~\eqref{eq:fx}, \eqref{eq:et}. 
Let $q$ be the ideal quadrilateral with sides $a,b,c,d$ and diagonals
$e,e'$. 
Then, equating the parallel transport operators along edges
of $q$, we
easily obtain eqns~\eqref{eq:ep}, \eqref{eq:ee}. For
example, by Lemma~\ref{l:11}, the parallel transport operator along
the short edge of 
$\GRF{\alpha,\tau'}$,
corresponding to the corner of $t_1'$ opposite to $e'$, is described by
the number 
\[
\EPS'(t_1')\frac{f'(e')}{f'(a)f'(d)}=\EPS'(t_1')\frac{f'(e')}{f(a)f(d)}
\]
On the other hand, the same path in $\GRF{\alpha,\tau}$ is represented
by composition of 
two short edges, corresponding to the corner of $t_1$ opposite to $b$
and the corner of $t_2$ opposite to $c$. Thus, the corresponding number
is given by the sum
\[
\EPS(t_1)\frac{f(b)}{f(a)f(e)}+\EPS(t_2)\frac{f(c)}{f(d)f(e)}
\]
The two numbers are equal by virtue of eqn~\eqref{eq:ee}.
\end{proof}
To complete the proof of part (ii), we note that
eqn~\eqref{eq:ep} implies that the quantity
$ k=\frac12\sum_{t\in\bar{\tau}}\EPS(t)$ is independent of
$\tau$ and thus the sets 
$\PBL_k(\tau)=\coprod_{\EPS\in\sigma^{-1}(k)}\PBL_{\EPS}(\tau)$
are such that $\PBL_k(\tau)\cap\PBL_l(\tau')=\emptyset$ if $k\ne l$
and any $\tau,\tau'\in\IDT$.

\subsection{Proof of part (iii)}
This is a direct consequence of 
Propositions~\ref{p:3} and \ref{p:4}, which in the two cases when all
signs, associated with ideal triangles, coincide give the formulas
identical to those of Penner for the decorated Teichm\"uller space.

\end{document}